\documentclass[12pt]{amsart}
\usepackage[utf8]{inputenc}
\usepackage{amsthm}
\usepackage{amsmath}
\usepackage{amssymb}
\usepackage{comment}
\usepackage[margin=1.2in]{geometry}

\usepackage[T1]{fontenc}
\usepackage{enumitem}
\usepackage{setspace}
\usepackage{microtype}
\usepackage{cite}
\usepackage{tikz}
\usetikzlibrary{decorations.pathreplacing}
\usepackage{mathtools}

\allowdisplaybreaks

\usepackage[colorlinks=true, pdfstartview=FitV, linkcolor=blue, citecolor=blue, urlcolor=blue]{hyperref}

\newtheorem{theorem}{Theorem}[section]

\newtheorem{lemma}[theorem]{Lemma}

\newtheorem{corollary}[theorem]{Corollary}

\providecommand{\ol}{\overline}

\providecommand{\C}{\mathbb C}

\providecommand{\N}{\mathbb N}

\providecommand{\R}{\mathbb R}
\providecommand{\Z}{\mathbb Z}

\title[$\ell^p$-Improving for the Moment Curve]{$\ell^p$-Improving for Discrete Fractional Integrals and Averages on the Moment Curve}

\author{Jacob Glidewell}
\address{Dept. of Mathematics \\
Indiana University \\
 Bloomington, IN, 47405-7106, USA}
\email{jaglide@iu.edu}

\begin{document}

\begin{abstract}
We prove the $\ell^p$-improving estimates in the full range for discrete fractional integral operators on the moment curve (except at the endpoint) and discrete averages on the moment curve (with $N^{\epsilon}$-losses). This follows via essentially sharp Diophantine solution counting estimates for inhomogeneous Vinogradov systems with an odd number of variables. 
\end{abstract}

\maketitle

\section{Introduction}
Fix $d\ge 2$. Consider the moment curve $$\gamma(n) = (n, n^2, \ldots, n^d).$$ We consider the following discrete averaging and fractional integral operators: for $f:\Z^d\to \C$ and $x\in \Z^d$, $$A_Nf(x)\coloneqq \frac{1}{N}\sum_{n=1}^N f(x-\gamma(n)),$$ and for $0<\lambda<1$, $$I_\lambda f(x)\coloneqq \sum_{n=1}^\infty \frac{f(x-\gamma(n))}{n^\lambda}.$$ We state our main theorem.

\begin{theorem}\label{thm:MainEll^pImproving}
If $0<\lambda<1$ and $(p,q)\in [1,\infty]^2$ satisfies $$\frac{1}{q}< \frac{1}{p}-\frac{2(1-\lambda)}{d(d+1)},\quad \frac{1}{p}>1-\lambda, \quad \frac{1}{q}<\lambda,$$ then for any $f:\Z^d\to \C$, $$\|I_\lambda f\|_{\ell^q(\Z^d)}\lesssim_{\lambda} \|f\|_{\ell^p(\Z^d)}.$$ 
\end{theorem}

Observe that the condition for $\lambda$ and $(p,q)$ in Theorem \ref{thm:MainEll^pImproving} is equivalent to $$0<1-\lambda<\min \left\{\frac{d(d+1)}{2}\left(\frac{1}{p}-\frac{1}{q}\right),\; \frac{1}{p},\; 1-\frac{1}{q}\right\}.$$ A simple transference principle between $A_N$ and $I_\lambda$ shows the next corollary follows immediately from Theorem \ref{thm:MainEll^pImproving}. This answers \cite[Conjecture 1]{DHV22} in the case of the moment curve.

\begin{corollary} \label{cor:DiscreteAvgs}
    If $1<p\le q<\infty$ and $\epsilon>0$, then for any $f:\Z^d\to \C$, $$\|A_Nf\|_{\ell^q(\Z^d)}\lesssim_\epsilon N^{\epsilon}[N^{-\frac{d(d+1)}{2}\left(\frac{1}{p}-\frac{1}{q}\right)}+N^{\frac{1}{q}-1}+N^{-\frac{1}{p}}]\|f\|_{\ell^p(\Z^d)}.$$
\end{corollary}

 Define for $r\in \N$ and $z\in \Z^d$, $$\mathcal{N}_r(N; z) = \#\{(n_1,\ldots, n_r)\in ([-N,N]\cap \Z)^{r} : \sum_{i=1}^r (-1)^{i+1}\gamma(n_i) = z\}.$$ Kim \cite{Kim15} showed that the full range (excluding endpoints) of $\ell^p(\Z^d)\to \ell^q(\Z^d)$ boundedness for $I_\lambda$ is implied by the following counting estimate. In particular, this is \cite[Eq. (2.4)]{Kim15}. Note that we are counting solutions to an inhomogeneous system with an \textbf{odd} number of variables $n_i$ since $d^2+d-1$ is always odd.

\begin{theorem}\label{thm:OddsolutionCounting}
     For any $z\in \Z^d$, $\mathcal{N}_{d^2+d-1}(N; z)\lessapprox N^{\frac{d(d+1)}{2}-1}$, with the implicit constant independent of $z$.
\end{theorem}

This is the sharpest bound to expect for general $z\in \Z^d$. Indeed, if $z=\gamma(n_0)$ for some $n_0\in [-N,N]$, then the diagonal solutions $$n_1=n_2,\quad n_3=n_4,\quad \ldots \quad n_{d^2+d-3}=n_{d^2+d-2}$$ and $n_{d^2+d-1}=n_0$ contribute $\sim N^{\frac{d^2+d-2}{2}}= N^{\frac{d(d+1)}{2}-1}$ solutions. Moreover, following the development of this work, it came to the author's attention that \cite[Theorem 1.2]{Wooley23c} gives the desired bound in Theorem \ref{thm:OddsolutionCounting} for $z=0$. 

As per \cite[Conjecture 1]{Kim15}, it is also expected $I_\lambda$ should be bounded on $\ell^p(\Z^d)\to \ell^q(\Z^d)$ when we have equality $$\frac{1}{p}=\frac{1}{q}-\frac{2(1-\lambda)}{d(d+1)}.$$ However, this would require our estimate in Theorem \ref{thm:OddsolutionCounting} without an $N^{\epsilon}$-loss, which would follow, via our proof, from VMVT without any $N^{\epsilon}$-loss.

The outline of this note is as follows. For the remainder of the introduction, we will review the previous work on $\ell^p$-improving inequalities for discrete averages and fractional integrals, and state an application toward discrete averages on polynomials $Q:\Z\to \Z$. In Section \ref{sec:Thm2impliesThm1}, we will show how Theorem \ref{thm:MainEll^pImproving} follows from Theorem \ref{thm:OddsolutionCounting}. In Section \ref{sec:OddSolnCount}, we will give a proof of Theorem \ref{thm:OddsolutionCounting}, first for $d=3$ to clarify the main ideas and then in general $d\ge 2$. Finally, in Section \ref{sec:TransPrinciple}, we give the short transference principle to prove Corollary \ref{cor:DiscreteAvgs}.

\subsection{Notation}

Throughout, we use $A\lesssim B$ to mean $A\le CB$ for some constant $C>0$, possibly depending on $d, p, q$ but otherwise universal, and $A\sim B$ to mean $A\lesssim B$ and $B\lesssim A$ are both satisfied. Further, we write $A\lessapprox B$ if for every $\epsilon>0$, we have $A\lesssim_\epsilon N^{\epsilon} B$, where the implied constant depends on $\epsilon$.

We use $\Z=\{\ldots, -2, -1, 0, 1, 2, \ldots\}$ and $\N = \{1, 2, 3, \ldots\}$. If $E\subset \Z^r$ for some $r\in \N$, we write $\#E$ for its cardinality. Also, for $1\le p<\infty$, $$\|f\|_{\ell^p(\Z^d)}=\left(\sum_{x\in \Z^d}|f(x)|^p\right)^{1/p},$$ and $$\|f\|_\infty = \sup_{x\in \Z^d}|f(x)|.$$ For $\theta\in \R$, we write $e(\theta)=e^{2\pi i \theta}$. 

\subsection{Previous Work}
Discrete fractional integrals were initially studied in \cite{SW00}, \cite{SW02}, \cite{Oberlin01}, \cite{IW06}. subsequent work was done by \cite{Pierce11}, \cite{Pierce12} and \cite{Kim15}, connecting discrete fractional integrals with Vinogradov's mean value theorem and related Diophantine solution counting estimates. In recent years, there has also been a large body of work on $\ell^p$-improving estimates for discrete averaging operators on curves \cite{HKLMY20}, \cite{HLY2021}, \cite{DHV22}, \cite{Hughes22};  surfaces \cite{DDL21}; or other arithmetic sets such as lattice points on the sphere \cite{H2020}, \cite{KL2020}, \cite{HY26}; or the primes \cite{Anderson2020}, \cite{HKLY2020}, \cite{LMR2022}, \cite{GLM2022}. 

Our work directly builds on the investigations of Kim \cite{Kim15} and Dendrinos--Hughes--Vitturi \cite{DHV22}. In \cite{Kim15}, Kim provides a mechanism that converts counting estimates, such as our Theorem \ref{thm:OddsolutionCounting}, into $\ell^p(\Z^d)\to \ell^q(\Z^d)$ bounds for discrete fractional integrals. In \cite{DHV22}, progress in the subcritical regime for discrete curve averages is made by counting solutions to inhomogeneous Diophantine equations with an even number of variables. Both \cite{Kim15} and \cite{DHV22} are built on discrete formulations of Christ's \textit{method of refinements} \cite{Christ98}.

Since the work of \cite{DHV22}, inhomogeneous Vinogradov systems with an even number of variables have been further studied \cite{Hughes22}, \cite{BH22}, \cite{Wooley23b}, \cite{Wooley23a}, \cite{Wooley23c}, \cite{Wooley25}, \cite{brandesWooley2026}. Moreover, inhomogeneous systems were crucial in the recent work \cite{wooley2026mainconjecturevinogradovsubsystem} on Vinogradov subsystems. Compared with our proof of Theorem \ref{thm:OddsolutionCounting}, counting solutions to systems with an even number of variables is a much harder problem. One difference is that the odd case can count `diagonal' solutions while the even case must avoid these. See the discussion after Theorem \ref{thm:OddsolutionCounting} above.

\subsection{Applications}

By the argument given in \cite[Theorem 1.9]{HKLMY20}, we get $\ell^p$-improving for polynomials in $\Z$, albeit for a suboptimal range. Let $Q$ be a polynomial of degree $d$ mapping $\Z\to \Z$. Define the averaging operator associated to $Q$: for $f:\Z\to \C$ and $x\in \Z$, $$A_N^Qf(x)= \frac{1}{N}\sum_{n=1}^N f(x+Q(n)).$$

\begin{corollary}\label{cor:PolynomialProjection}
    For $2-\frac{2}{d(d+1)}\le p\le 2$, we have $$\|A_N^Qf\|_{\ell^{p^{\prime}}(\Z)}\lessapprox N^{-d\left(\frac{2}{p}-1\right)}\|f\|_{\ell^p(\Z)},$$ for any $f\in \ell^p(\Z)$, where the implicit constant may depend on $Q$.
\end{corollary}

This is not the optimal range of $p$ but improves the range in \cite{HKLMY20}. The endpoint exponent is expected to be $p=2-\frac{1}{d}$. 

Consider the cubes $Q(n)=n^3$. By the same arguments for Theorem \ref{thm:MainEll^pImproving}, the optimal $\ell^p$-improving range for $A_N^Q$ would follow by showing for a given $z\in \Z$, there are $\lessapprox N^2$ solutions $(n_1,\ldots, n_5)\in [-N,N]^5$ to the Diophantine equation $$z=n_1^3-n_2^3+n_3^3-n_4^3+n_5^3.$$

\subsection{Declarations}

\medskip
\noindent
{\it Use of AI}. Artificial Intelligence was used at the end of the work for aid in literature review. This note was solely written and edited by the author. All mathematical claims are the work and responsibility of the author.

\medskip
\noindent
{\it Acknowledgments}. The author would like to thank Kevin Hughes for feedback on an earlier draft of this note and education on the previous work related to discrete $\ell^p$-improving estimates and inhomogeneous Vinogradov systems.

\section{Proof of Theorem \ref{thm:MainEll^pImproving} from Theorem \ref{thm:OddsolutionCounting}}\label{sec:Thm2impliesThm1}

The goal of this section is to show Theorem \ref{thm:MainEll^pImproving} follows from Theorem \ref{thm:OddsolutionCounting}. We will apply \cite[Theorem 1.1]{Kim15} for a range of $\lambda$ and then use H\"{o}lder's inequality for the remaining range of $\lambda$. We restate \cite[Theorem 1.1]{Kim15} for the sake of the reader.

\begin{theorem}\label{thm:Kim'sTheorem}
    \cite{Kim15} Suppose for some $r\in \N$ and $\delta>0$, $$\sup_{z\in \Z^d} \mathcal{N}_r(N;z)\lessapprox N^{r-\delta}.$$ Let $s\in \N$ such that either $r=2s$ or $r=2s-1$. Then, $I_\lambda$ is bounded $\ell^p(\Z^d)\to \ell^q(\Z^d)$ if $1-\frac{\delta}{r}<\lambda<1$ and $(p,q)\in [1,\infty]^2$ satisfy $$\frac{1}{q}< \frac{1}{p}-\frac{(1-\lambda)}{\delta},\quad \frac{1}{p}>\frac{s}{\delta}(1-\lambda), \quad \frac{1}{q}<1-\frac{s}{\delta}(1-\lambda).$$
\end{theorem}

Hence, given Theorem \ref{thm:OddsolutionCounting}, we may apply Theorem \ref{thm:Kim'sTheorem} with $r=d^2+d-1$ and $\delta=s=\frac{d(d+1)}{2}$ to get Theorem \ref{thm:MainEll^pImproving} in the range $$\lambda_d\coloneqq 1-\frac{d(d+1)}{2(d^2+d-1)}<\lambda<1.$$ As observed by Kim \cite{Kim15}, the complementary range $0<\lambda \le \lambda_d$ follows by interpolating the range $\lambda_d < \lambda < 1$ with the trivial estimate $$\|I_\lambda f\|_\infty \le \|f\|_{\ell^1(\Z^d)}.$$ We give a short proof of this fact via H\"{o}lder's inequality. Fix $0<\lambda \le \lambda_d$ and $(p,q)\in [1,\infty]^2$ that satisfies $$\frac{1}{q}< \frac{1}{p}-\frac{2(1-\lambda)}{d(d+1)},\quad \frac{1}{p}>1-\lambda, \quad \frac{1}{q}<\lambda.$$ Pick $\mu>\lambda_d$ but sufficiently close to $\lambda_d$ and $\theta = \frac{\lambda}{\mu}\in (0,1)$. Define $p_1$ and $q_1$ by $$\frac{1}{p}= (1-\theta)+\frac{\theta}{p_1}, \quad \frac{1}{q}= \frac{\theta}{q_1}.$$ Notice that $(p_1,q_1)$ satisfies $$\frac{1}{q_1}< \frac{1}{p_1}-\frac{2(1-\mu)}{d(d+1)},\quad \frac{1}{p_1}>1-\mu, \quad \frac{1}{q_1}<\mu.$$ Indeed, $$\frac{1}{p_1}-1+\mu = \frac{1}{\theta}\left(\frac{1}{p}-1+\lambda\right)>0,$$ $$\frac{1}{q_1}-\mu =\frac{1}{\theta}\left(\frac{1}{q}-\lambda\right)<0,$$ and \begin{align*}
    \frac{1}{p_1}-\frac{1}{q_1}&= \frac{1}{\theta}\left(\frac{1}{p}-\frac{1}{q}-1+\theta\right)=\frac{1}{\theta}\left(\frac{1}{p}-\frac{1}{q}-1+2\lambda\right)+1-2\mu\\
    &= \frac{1}{\theta}\left(\frac{1}{p}-\frac{1}{q}-1+2\lambda\right)+\frac{2(1-\mu)}{d(d+1)}-\frac{2(d^2+d-1)}{d(d+1)}(\mu-\lambda_d)\\
    &>\frac{2(1-\mu)}{d(d+1)},
\end{align*} for $\mu$ sufficiently close to $\lambda_d$. This gives that $I_\mu$ is bounded $\ell^{p_1}(\Z^d)\to \ell^{q_1}(\Z^d)$. 

Fix $f\in \ell^p(\Z^d)$. Without loss of generality, we may assume $f\ge 0$ by replacing $f$ with $|f|$. First, observe the pointwise inequality \begin{equation}\label{eq:pwHolder}
    I_\lambda f(x) \le (I_\mu [f^{p/p_1}](x))^{\theta}\|f\|_{\ell^p(\Z^d)}^{p(1-\theta)}.
\end{equation} Indeed, by H\"{o}lder's inequality, \begin{align*}
    I_\lambda f(x) &= \sum_{n\ge 1} \frac{f(x-\gamma(n))}{n^{\lambda}} = \sum_{n\ge 1}\frac{f(x-\gamma(n))^{\theta\frac{p}{p_1}+p(1-\theta)}}{n^{\lambda}}\\
    &\le \left(\sum_{n\ge 1} \frac{f(x-\gamma(n))^{p/p_1}}{n^{\mu}}\right)^{\theta}\left(\sum_{n\ge 1} f(x-\gamma(n))^{p}\right)^{1-\theta}.
\end{align*}
Now, from our pointwise inequality \eqref{eq:pwHolder}, we may complete the proof. We get $$\|I_\lambda f\|_{\ell^q(\Z^d)}\le \|f\|_{\ell^p(\Z^d)}^{p(1-\theta)}\|I_\mu[f^{p/p_1}]\|_{\ell^{q\theta}(\Z^d)}^{\theta}= \|f\|_{\ell^p(\Z^d)}^{p(1-\theta)}\|I_\mu[f^{p/p_1}]\|_{\ell^{q_1}(\Z^d)}^{\theta}.$$ Since $I_\mu$ is bounded $\ell^{p_1}(\Z^d)\to \ell^{q_1}(\Z^d)$, we get $$\|I_\mu[f^{p/p_1}]\|_{\ell^{q_1}(\Z^d)}^{\theta}\lesssim_\mu \|f^{p/p_1}\|_{\ell^{p_1}(\Z^d)}^{\theta}= \|f\|_{\ell^p(\Z^d)}^{p\theta/p_1}.$$ Hence, we are finished with $$\|I_\lambda f\|_{\ell^q(\Z^d)}\lesssim \|f\|_{\ell^p(\Z^d)}^{p\theta/p_1}\|f\|_{\ell^p(\Z^d)}^{p(1-\theta)}=\|f\|_{\ell^{p}(\Z^d)}.$$

\section{Proof of Theorem \ref{thm:OddsolutionCounting}}\label{sec:OddSolnCount}
In general, our proof is similar to that in \cite{BH22}, \cite{Wooley23c}, \cite{wooley2026mainconjecturevinogradovsubsystem}: we use a shifting/averaging argument followed by two applications of VMVT of Wooley \cite{Wooley_2016} and Bourgain--Demeter--Guth \cite{BDG2016}. We summarize this estimate as follows. For $t\in [0,1]^d$, let $$S_N(t)= \sum_{|n|\le N}e(t\cdot \gamma(n)).$$

\begin{theorem}\label{thm:VMVT}
    (\cite{Wooley_2016} and \cite{BDG2016}) $$\int_{\mathbb{T}^d} |S_N(t)|^{d(d+1)}\; dt \lessapprox N^{\frac{d(d+1)}{2}}.$$
\end{theorem}

\subsection{Proof for $d=3$.}

For the sake of clarity, we first present the argument for the moment curve with $d=3$: $$\gamma(n)= (n, n^2, n^3).$$ It suffices to see that for any $z\in \Z^3$, $$\mathcal{N}_{11}(N;z)\lessapprox N^{5}.$$ 

Now, fix $z=(z_1,z_2, z_3)\in \Z^3$ and $h\in [-N,N]\cap \Z$. Consider a solution $$(n_1, \ldots, n_{11})$$ which contributes to $\mathcal{N}_{11}(N;z)$ and consider the shift $m_i=n_i+h$. Then, $$(m_1, \ldots, m_{11})\in [-2N, 2N]^{11}$$ contributes to $\mathcal{N}_{11}(2N; T_h(z))$ where $$T_h(z)= (z_1+h, z_2+2hz_1+h^2, z_3+3hz_2+3h^2z_1+h^3).$$ Averaging over $h$ gives $$\mathcal{N}_{11}(N; z)\le \frac{1}{N}\sum_{|h|\le N}\mathcal{N}_{11}(2N; T_h(z)).$$ Interpreting this count via mean values gives \begin{align*}
        \sum_{|h|\le N}\mathcal{N}_{11}(2N; T_h(z)) &= \sum_{|h|\le N} \int_{\mathbb{T}^{3}} (S_{2N}(t))^{6}(\ol{S_{2N}(t)})^{5}e(-T_h(z)\cdot t)\; dt.
    \end{align*}
    By H\"{o}lder's inequality, we bound this by $$\left(\int_{\mathbb{T}^3}|S_{2N}(t)|^{12}\; dt\right)^{\frac{11}{12}}\left(\int_{\mathbb{T}^3}\left|\sum_{h}e(T_h(z)\cdot t)\right|^{12}\; dt\right)^{\frac{1}{12}}.$$ Using Theorem \ref{thm:VMVT}, the first term is bounded by $$\left(\int_{\mathbb{T}^3}|S_{2N}(t)|^{12}\; dt\right)^{\frac{11}{12}}\lessapprox N^{\frac{11}{2}}.$$ By a change-of-variables with Jacobian $1$, we may apply Theorem \ref{thm:VMVT} to the second term as well! Indeed, \begin{align*}
        T_h(z)\cdot t &= (z_1t_1+z_2t_2+z_3t_3)\\
        &\quad +h(t_1+2z_1t_2+3z_2t_3)\\
        &\quad +h^2(t_2+3z_1t_3)\\
        &\quad +h^3t_3.
    \end{align*}This gives $$\int_{\mathbb{T}^3}\left|\sum_{h}e(T_h(z)\cdot t)\right|^{12}\; dt = \int_{\mathbb{T}^3}\left|\sum_{h}e((h, h^2, h^3)\cdot At)\right|^{12}\; dt,$$ where $$A= \begin{bmatrix}
        1 & 2z_1 & 3z_2\\
        0 & 1 & 3z_1\\
        0 & 0 & 1
    \end{bmatrix}.$$ Since $\det (A)=1$ and in particular, $A\in \mathrm{SL}_3(\Z)$, we let $u=At$. This becomes $$\left(\int_{\mathbb{T}^3}\left|\sum_{|h|\le N}e(\gamma(h)\cdot u)\right|^{12}\; du\right)^{\frac{1}{12}}\lessapprox N^{1/2}.$$ Hence, we combine the previous estimates to get $$\mathcal{N}_{11}(N; z)\lessapprox \frac{1}{N}\cdot  N^{\frac{11}{2}}\cdot N^{1/2}= N^{5}.$$

\subsection{Proof for general $d\ge 2$.}
The case $d=2$ is classical \cite{Oberlin01}; however, it fits into our larger framework. We are considering the curve $$\gamma(n)=(n, n^2,\ldots, n^d).$$ It suffices to see that for any $z\in \Z^d$, $$\mathcal{N}_{d^2+d-1}(N; z)\lessapprox N^{\frac{d(d+1)}{2}-1}=N^{\frac{d^2+d-2}{2}}.$$

Now, fix $z=(z_1,\ldots, z_d)\in \Z^d$ and $h\in [-N,N]\cap \Z$.

\begin{lemma}
    Consider a solution $$(n_1, \ldots, n_{d^2+d-1})\in ([-N,N]\cap \Z)^{d^2+d-1}$$ which contributes to $\mathcal{N}_{d^2+d-1}(N; z)$ and consider the shift $m_i=n_i+h$. Then, $$(m_1, \ldots, m_{d^2+d-1})\in ([-2N, 2N]\cap\Z)^{d^2+d-1}$$ contributes to $\mathcal{N}_{d^2+d-1}(2N; T_h(z))$ where $$T_h(z)= (z_1+h, z_2+2hz_1+h^2, z_3+3hz_2+3h^2z_1+h^3, \ldots)\in \Z^d.$$ That is, $[T_h(z)]_k= h^k+\sum_{j=0}^{k-1}\binom{k}{j}h^jz_{k-j}$.
\end{lemma}

\begin{proof}
    We need to show that for each $1\le k\le d$, \begin{equation*}\label{eq:inductonk}
        \sum_{i=1}^{d^2+d-1}(-1)^{i+1}(n_i+h)^{k}= [T_h(z)]_k.
    \end{equation*} We will use the convention $z_0=1$. Then, \begin{align*}
        \sum_{i=1}^{d^2+d-1}(-1)^{i+1}(n_i+h)^{k}&= \sum_{i=1}^{d^2+d-1}(-1)^{i+1}\sum_{j=0}^{k}\binom{k}{j}n_i^{k-j}h^j\\
        &=\sum_{j=0}^k \binom{k}{j}h^j\sum_{i=1}^{d^2+d-1}(-1)^{i+1}n_i^{k-j}\\
        &= \sum_{j=0}^{k} \binom{k}{j}h^j z_{k-j} = [T_h(z)]_k.
    \end{align*}
\end{proof}
 
 Averaging over $h$ gives $$\mathcal{N}_{d^2+d-1}(N; z)\le \frac{1}{N}\sum_{|h|\le N}\mathcal{N}_{d^2+d-1}(2N; T_h(z)).$$ Interpreting this count via mean values gives \begin{align*}
        \sum_{|h|\le N}\mathcal{N}_{d^2+d-1}(2N; T_h(z)) &= \sum_{|h|\le N} \int_{\mathbb{T}^{d}} (S_{2N}(t))^{\frac{d(d+1)}{2}}(\ol{S_{2N}(t)})^{\frac{d(d+1)}{2}-1}e(-T_h(z)\cdot t)\; dt.
    \end{align*} By H\"{o}lder's inequality, we bound this by $$\left(\int_{\mathbb{T}^d}|S_{2N}(t)|^{d(d+1)}\; dt\right)^{\frac{d^2+d-1}{d^2+d}}\left(\int_{\mathbb{T}^d}\left|\sum_{h}e(T_h(z)\cdot t)\right|^{d(d+1)}\; dt\right)^{\frac{1}{d^2+d}}.$$ Using Theorem \ref{thm:VMVT}, the first term is bounded by $$\left(\int_{\mathbb{T}^d}|S_{2N}(t)|^{d(d+1)}\; dt\right)^{\frac{d^2+d-1}{d^2+d}}\lessapprox N^{\frac{d^2+d-1}{2}}.$$ By a change-of-variables with Jacobian $1$, we may apply Theorem \ref{thm:VMVT} to the second term as well! Indeed, by the same computation as in the $d=3$ case, $$T_h(z)\cdot t= z\cdot t + \gamma(h)\cdot At,$$ which gives $$\int_{\mathbb{T}^d}\left|\sum_{h}e(T_h(z)\cdot t)\right|^{d(d+1)}\; dt = \int_{\mathbb{T}^d}\left|\sum_{h}e(\gamma(h)\cdot At)\right|^{d(d+1)}\; dt,$$ where $$A= \begin{bmatrix}
        1 & 2z_1 & 3z_2 & \cdots & dz_{d-1}\\
        0 & 1 & 3z_1 & \cdots & \binom{d}{2}z_{d-2}\\
        0 & 0 & 1 & \cdots & \binom{d}{3}z_{d-3}\\
        \vdots & \vdots & \vdots & \cdots & \vdots\\
        0 & 0 & 0 &\cdots & 1
    \end{bmatrix}.$$ Since $A\in \mathrm{SL}_d(\Z)$, we let $u=At$.  This becomes $$\left(\int_{\mathbb{T}^d}\left|\sum_{|h|\le N}e(\gamma(h)\cdot u)\right|^{d(d+1)}\; du\right)^{\frac{1}{d^2+d}}\lessapprox N^{1/2}.$$ Hence, we combine the previous estimates to get $$\mathcal{N}_{d^2+d-1}(N; z)\lessapprox \frac{1}{N}\cdot  N^{\frac{d^2+d-1}{2}}\cdot N^{1/2}= N^{\frac{d^2+d-2}{2}}.$$

\section{Proof of Corollary \ref{cor:DiscreteAvgs}}\label{sec:TransPrinciple}

We follow the same transference principle as in \cite{HKLMY20}. Note that for $f:\Z^d\to \C$ and $x\in \Z^d$, \begin{equation}\label{eq:pointwiseDom}
    |A_Nf(x)|\le \frac{1}{N}\sum_{n=1}^{N}|f|(x-\gamma(n))\le \frac{1}{N^{1-\lambda}}\sum_{n=1}^N \frac{|f|(x-\gamma(n))}{n^\lambda} \le \frac{1}{N^{1-\lambda}}I_\lambda|f|(x).
\end{equation} Fix $1<p\le q<\infty$ and $\epsilon>0$ sufficiently small. If $p=q$, then Corollary \ref{cor:DiscreteAvgs} is trivial. Assume $p<q$. By \eqref{eq:pointwiseDom}, we have $$\|A_N f\|_{\ell^q(\Z^d)}\le N^{-(1-\lambda)}\|I_\lambda|f|\|_{\ell^q(\Z^d)}.$$ Now, pick $\lambda\in (0,1)$ such that $$1-\lambda=-\epsilon+ \min \left\{\frac{d(d+1)}{2}\left(\frac{1}{p}-\frac{1}{q}\right), \frac{1}{p}, 1-\frac{1}{q}\right\}.$$ Now, we use Theorem \ref{thm:MainEll^pImproving} to get $$\|A_N f\|_{\ell^q(\Z^d)}\lesssim_\epsilon N^{\epsilon}[N^{-\frac{d(d+1)}{2}\left(\frac{1}{p}-\frac{1}{q}\right)}+N^{\frac{1}{q}-1}+N^{-\frac{1}{p}}]\|f\|_{\ell^p(\Z^d)}.$$
\bibliographystyle{plain}
\bibliography{refImproving}

\end{document}